\documentclass[epsfig,latexsym,amsfonts,twoside]{article}
\usepackage{amssymb,amsmath,amscd}

\def\part#1{\frac{\partial\phantom{#1}}{\partial#1}}
\newtheorem{thm}{Theorem}
\newtheorem{prp}[thm]{Proposition}
\newtheorem{lem}[thm]{Lemma}
\newtheorem{cor}[thm]{Corollary}

\newenvironment{prf}{\begin{trivlist}\item[]{\bf Proof} }%
{\hfill $\Box$ \end{trivlist}}
\newenvironment{dfn}{\begin{trivlist}\item[]{\bf Definition}\em }%
{\end{trivlist}}
\newenvironment{rmk}{\begin{trivlist}\item[]{\bf Remark} }%
{\end{trivlist}}
\newenvironment{exm}{\begin{trivlist}\item[]{\bf Example} }%
{\end{trivlist}}

\def\Z{\ifmmode{{\mathbb Z}}\else{${\mathbb Z}$}\fi}
\def\Q{\ifmmode{{\mathbb Q}}\else{${\mathbb Q}$}\fi}
\def\C{\ifmmode{{\mathbb C}}\else{${\mathbb C}$}\fi} 
\def\P{\ifmmode{{\mathbb P}}\else{${\mathbb P}$}\fi} 

\def\H{\ifmmode{{\mathrm H}}\else{${\mathrm H}$}\fi} 

\def\B{\ifmmode{{\cal B}}\else{${\cal B}$}\fi} 
\def\E{\ifmmode{{\cal E}}\else{${\cal E}$}\fi} 
\def\F{\ifmmode{{\cal F}}\else{${\cal F}$}\fi} 
\def\I{\ifmmode{{\cal I}}\else{${\cal I}$}\fi} 
\def\K{\ifmmode{{\cal K}}\else{${\cal K}$}\fi} 
\def\L{\ifmmode{{\cal L}}\else{${\cal L}$}\fi} 
\def\M{\ifmmode{{\cal M}}\else{${\cal M}$}\fi} 
\def\N{\ifmmode{{\cal N}}\else{${\cal N}$}\fi} 
\def\O{\ifmmode{{\cal O}}\else{${\cal O}$}\fi} 
\def\U{\ifmmode{{\cal U}}\else{${\cal U}$}\fi}
\def\X{\ifmmode{{\cal X}}\else{${\cal X}$}\fi} 

\def\Br{\ifmmode{{\mathrm{Br}}}\else{${\mathrm{Br}}$}\fi} 
\def\OG{\ifmmode{\widetilde{\cal M}_4}\else{$\widetilde{\cal M}_4$}\fi} 
\def\D{\ifmmode{{\cal D}_{\mathrm{coh}}^b}\else{${{\cal
    D}_{\mathrm{coh}}^b}$}\fi}
\def\Shah{\ifmmode{\amalg\hspace*{-3.5pt}\amalg}\else{$\amalg\hspace*{-3.5pt}\amalg$}\fi}

\begin{document}

\title{Fibrations on four-folds with trivial canonical
  bundles\footnote{2000 {\em Mathematics Subject Classification.\/}
    14J35; 14D06, 14J32, 53C26.}}
\author{Justin Sawon}
\date{April, 2009}
\maketitle

\begin{abstract}
Four-folds with trivial canonical bundles are divided into six classes
according to their holonomy group. We consider examples that are
fibred by abelian surfaces over the projective plane. We construct such
fibrations in five of the six classes, and prove that there is no
such fibration in the sixth class. We classify all such fibrations
whose generic fibre is the Jacobian of a genus two curve.
\end{abstract}

\section{Introduction}

A holomorphic symplectic four-fold $X$ is a compact K{\"a}hler
four-fold admitting a holomorphic two-form
$\sigma\in\mathrm{H}^0(X,\Omega^2)$ which is non-degenerate, i.e.,
$\sigma^{\wedge 2}$ trivializes the canonical bundle
$K_X=\Omega^4$. The Calabi-Yau metric on $X$ has holonomy group
contained in $\mathrm{Sp}(2)$. We say that $X$ is irreducible if it is
simply-connected and $\sigma$ is unique up to scale; in this case the
holonomy group is equal to
$\mathrm{Sp}(2)$. Matsushita~\cite{matsushita99,matsushita00i,matsushita00ii}
showed that if a projective irreducible holomorphic symplectic
four-fold $X$ admits a fibration then the generic fibre must be an
abelian surface and the base must be isomorphic to $\P^2$ (see also
Markushevich~\cite{markushevich96} for the statement about the base).

In the opposite direction, one might ask: ``when is a fibration
$X\rightarrow\P^2$ by abelian surfaces an irreducible holomorphic
symplectic four-fold?'' (We mean here a fibration whose generic fibre
is an abelian surface; we always allow singular fibres.) A necessary
condition is that the canonical bundle $K_X$ of $X$ must be
trivial. Now there is a formula to compute $K_X$ due to Fujino and
Mori~\cite{fm00} which
generalizes Kodaira's formula for the canonical bundle of an elliptic
surface; this could be used to check whether $K_X$ is
trivial. However, this is not enough: as we show in this article,
there are numerous examples of four-folds with trivial canonical
bundles admitting fibrations by abelian surfaces over $\P^2$ which are
{\em not\/} irreducible holomorphic symplectic four-folds.

Four-folds with trivial canonical bundles admit Calabi-Yau metrics;
they can then be divided into six classes according to the restricted
holonomy group of this metric. Irreducible holomorphic symplectic
four-folds have holonomy $\mathrm{Sp}(2)$; four of the other classes
also contain examples of fibrations by abelian surfaces over $\P^2$.
\begin{thm}
\begin{enumerate}
\item Let $G$ be the trivial group, $\mathrm{SU}(2)$,
  $\mathrm{SU}(3)$, $\mathrm{SU}(4)$, or $\mathrm{Sp}(2)$. Then there
  exist four-folds $X$ with $K_X$ trivial and holonomy $G$, that are
  fibred by abelian surfaces over $\P^2$.
\item There does not exist a four-fold $X$ with $K_X$ trivial and
  holonomy $\mathrm{SU}(2)\times\mathrm{SU}(2)$, that is fibred by
  abelian surfaces over $\P^2$.
\end{enumerate}
\end{thm}
The six different classes are treated in Sections~3 to~8; part two of
the theorem is proved in Section~5 (Corollary~\ref{impossible1} and
Proposition~\ref{impossible2}) and the examples are constructed in the
other sections.

In addition, suppose that $X$ is fibred by Jacobians of genus two
curves (roughly this means that the fibres are principally polarized
abelian surfaces). Following ideas of
Markushevich~\cite{markushevich95}, we classify such fibrations. 
\begin{thm}
Let $\mathcal{C}\rightarrow\P^2$ be a family of genus two curves with
`mild degenerations' (defined in Section~9), and suppose that the
compactified relative Jacobian
$$X:=\overline{\mathrm{Jac}}^0(X/\P^2)$$
has trivial canonical bundle. Then $X$ is either
\begin{enumerate}
\item a Calabi-Yau four-fold, with holonomy $\mathrm{SU}(4)$,
  contained in a connected $75$-dimensional family, 
\item or an irreducible holomorphic symplectic four-fold, with
  holonomy $\mathrm{Sp}(2)$, contained in a connected $19$-dimensional
  family.
\end{enumerate}
\end{thm}
The proof appears in Section~9 (Propositions~\ref{no0-3},~\ref{CY4-fold},
and~\ref{hilb2K3}). Markushevich~\cite{markushevich95}
proved that the holomorphic symplectic four-folds arising in case two
are precisely the Beauville-Mukai integrable
systems~\cite{beauville99}, i.e., the family
of curves $\mathcal{C}\rightarrow\P^2$ must be a complete linear system
of genus two curves on a K3 surface. We include a slight
simplification of his argument here for completeness. Markushevich's
theorem has been generalized to fibrations on higher-dimensional
holomorphic symplectic manifolds by the author~\cite{sawon08iv}.

In~\cite{sawon04, sawon08iii} the author described twisted
Fourier-Mukai transforms between Lagrangian fibrations
$X\rightarrow\P^n$ on holomorphic symplectic manifolds and their dual
fibrations $\widehat{X}\rightarrow\P^n$. The twisting is by a gerbe,
which lies in the analytic Brauer group $\H^2(\widehat{X},\O^*)$. For
holomorphic symplectic manifolds, this group is connected. Similarly,
one could construct twisted Fourier-Mukai transforms between the
fibrations constructed in this article and their duals; this may
reveal new phenomena, since the Brauer group is discrete in some
cases.

The author would like to thank Rick Miranda and Jason Starr for useful
conversations, and Markushevich whose papers~\cite{markushevich95,
  markushevich96} provided the basis for much of this work
(particularly Section~9). The author is grateful for the hospitality
of the Korea Institute for Advanced Studies where parts of this
article were completed.

\section{Four-folds with trivial canonical bundles}

\subsection{Holonomy classification}

We will study fibrations over the projective plane whose
generic fibres are abelian surfaces and whose total spaces have
trivial canonical bundles.

\begin{dfn}
Throughout this article $X$ will denote a compact K{\"a}hler manifold of
dimension four, with trivial canonical bundle $K_X\cong\O_X$. We will
assume
\begin{enumerate}
\item there is a surjective morphism $\pi:X\rightarrow\P^2$ with connected
  fibres,
\item the generic fibre $X_t:=\pi^{-1}(t)$ is an abelian surface,
\item $\pi:X\rightarrow\P^2$ is a flat fibration, and in particular
  equidimensional,
\item there exists a smooth section $s:\P^2\rightarrow X$ of $\pi$.
\end{enumerate}
\end{dfn}

\begin{rmk}
It would be desirable to relax the last two assumptions. In principal,
given a fibration $X\rightarrow\P^2$ that {\em does not\/} admit a
section, there should be an associated fibration $X^0\rightarrow\P^2$ that
does admit a section such that $X$ and $X^0$ are locally isomorphic as
fibrations. Over smooth fibres, $X^0$ is the relative Albanese of
$X/\P^2$, or equivalently the ``double-dual'', i.e., relative Picard
scheme of the relative Picard scheme of $X/\P^2$. In practice, the
existence of singular fibres makes the construction of $X^0$ from $X$
difficult (though not impossible: see~\cite{sawon04,sawon08iii}).
\end{rmk}

\begin{lem}
If $X$ satisfies the above assumptions then $X$ is projective.
\end{lem}

\begin{prf}
When $X$ is a holomorphic symplectic four-fold, this follows from an 
argument of Campana and Oguiso (Proposition~3.2 in~\cite{oguiso07}; 
see also Lemma~2 in~\cite{sawon09}). The same argument works in this
more general setting.
\end{prf}

Since $X$ is compact and K{\"a}hler, with trivial canonical bundle,
Yau's Theorem implies that $X$ admits a Ricci-flat metric
$g$. The restricted holonomy $\mathrm{Hol}^0(g)$ is the group of
holonomy maps obtained by parallel transport around null-homotopic
paths. Since $g$ is Ricci-flat, $\mathrm{Hol}^0(g)$ must be contained
in $\mathrm{SU}(4)$ and Berger's classification of holonomy groups
(see Theorem~3.4.1 of Joyce~\cite{joyce00}) yields the following six
classes.

\begin{prp}
\label{berger}
The restricted holonomy $\mathrm{Hol}^0(g)$ of the Ricci-flat metric
on $X$ must be one of the following:
\begin{enumerate}
\item the trivial group; there is a finite unramified cover
  $\tilde{X}$ of $X$ such that $\tilde{X}$ is an abelian four-fold,
\item $\mathrm{SU}(2)$; there is a finite unramified cover $\tilde{X}$
  of $X$ such that $\tilde{X}=A\times S$, where $A$ and $S$ are
  abelian and K3 surfaces respectively,
\item $\mathrm{SU}(2)\times\mathrm{SU}(2)$; either $X$ itself or an
  unramified double cover $\tilde{X}$ of $X$ is the product
  $S_1\times S_2$ of two K3 surfaces,
\item $\mathrm{SU}(3)$; there is a finite unramified cover $\tilde{X}$
  of $X$ such that $\tilde{X}=E\times Y$, where $E$ is an elliptic
  curve and $Y$ is a Calabi-Yau three-fold,
\item $\mathrm{SU}(4)$; $X$ is a Calabi-Yau
  four-fold,
\item $\mathrm{Sp}(2)$; $X$ is an irreducible
  holomorphic symplectic four-fold.
\end{enumerate}
\end{prp}

\begin{prf}
As mentioned above, the possibilities for the restricted holonomy
group follow from Berger's classification. Moreover, the Bogomolov
Decomposition Theorem~\cite{bogomolov74} (see also Proposition~6.2.2
of Joyce~\cite{joyce00}) states that $X$ admits a finite unramified
cover $\tilde{X}$ which is a product manifold
$$\tilde{X}= T\times Y_1\times\ldots\times Y_p\times
Z_1\times\ldots\times Z_q.$$
In this decomposition, $T$ is a flat complex torus (necessarily an
abelian variety since $X$ is projective), the $Y_i$s are Calabi-Yau
manifolds with holonomy $\mathrm{SU}(n_i)$, and the $Z_i$s are
irreducible holomorphic symplectic manifolds with holonomy
$\mathrm{Sp}(m_i)$. The small dimensions lead to the six possibilities
listed above.

To prove $\tilde{X}=X$ in classes 5 and 6, and either $\tilde{X}=X$
or $\tilde{X}$ is a double cover of $X$ in class 3, we compare
$\chi(\O_X)$ to $\chi(\O_{\tilde{X}})$. If $\tilde{X}\rightarrow X$ is
an unramified $d$-to-1 cover then
$$\chi(\O_{\tilde{X}})=d\chi(\O_X).$$
Moreover, any holomorphic $k$-form on $X$ will lift to a holomorphic
$k$-form on $\tilde{X}$, so
$$h^k(\O_{\tilde{X}})=h^0(\Omega^k_{\tilde{X}})\geq
h^0(\Omega^k_X)=h^k(\O_X)$$
for all $k$. We also have $h^4(\O_X)=h^0(\Omega^4_X)=1$ since
$\Omega^4_X=K_X$ is trivial. One easily concludes that $d$ must equal
one in classes 5 and 6, and must equal one or two in class 3.
\end{prf}

\begin{exm}
We describe a four-fold $X$ with trivial canonical bundle which
has $S_1\times S_2$ as an unramified double cover, as in class 3. Let
$S_i$ be a double cover of an Enriques surface $T_i$, with covering
involution $\gamma_i$. Then $(\gamma_1,\gamma_2)$ is a fixed-point
free involution on $S_1\times S_2$, and we let 
$$\tilde{X}:=S_1\times S_2\rightarrow X:=S_1\times
S_2/\langle(\gamma_1,\gamma_2)\rangle.$$
Denote by $\sigma_i$ the holomorphic symplectic form on $S_i$, and by
$p_i^*\sigma_i$ the pull-back to $\tilde{X}=S_1\times S_2$, where
$p_i$ is projection to the $i$th factor. Then $\gamma_i$ acts on
$\sigma_i$ by $-1$, so $(\gamma_1,\gamma_2)$ acts on $p_i^*\sigma_i$
by $-1$. Neither two-form $p_i^*\sigma_i$ descends to $X$, but their
product $p_1^*\sigma_1\wedge p_2^*\sigma_2$ does descend to give a
trivialization of $K_X$.

We conclude that
$$h^k(\O_X)=\left\{ \begin{array}{ll}
  1 & \mbox{if }k=0\mbox{ or }4, \\
  0 & \mbox{otherwise.} \\
\end{array}\right.$$
If the holonomy group $\mathrm{Hol}(g)$ of $X$
were the same as the restricted holonomy group
$\mathrm{Hol}^0(g)\cong\mathrm{SU}(2)\times\mathrm{SU}(2)$ then there
would be parallel (and hence holomorphic) two-forms on $X$; but
$h^0(\Omega^2_X)=h^2(\O_X)=0$. Therefore $\mathrm{Hol}(g)$ is strictly
bigger than $\mathrm{Hol}^0(g)$. On the other hand, there is a
surjective group homomorphism
$$\pi_1(X)\cong\Z/2\Z \rightarrow \mathrm{Hol}(g)/\mathrm{Hol}^0(g)$$
(see Proposition~2.2.6 of Joyce~\cite{joyce00}). An example of such a
group is
$$\mathrm{Hol}(g)=\left\{\left. \left(\begin{array}{cc}
  A & 0 \\
  0 & B \\
\end{array}\right) \right |
  A,B\in\mathrm{U}(2),\mathrm{det}A=\mathrm{det}B=\pm 1
  \right\}\subset\mathrm{SU}(4).$$
\end{exm}

\begin{rmk}
It is enough that just one of the involutions $\gamma_i$ be
fixed-point free, for then $(\gamma_1,\gamma_2)$ will be fixed-point
free. For example, $S_1$ could be the unramified double cover of an
Enriques surface $T_1$ and $S_2$ could be the double cover of $\P^2$
ramified over a sextic.
\end{rmk}

This construction exhibits all the possibilities in class 3, as
we now show.

\begin{prp}
\label{Enriques}
Let $X$ be a four-fold with trivial canonical bundle which has the
product $S_1\times S_2$ of two K3 surfaces as an unramified double
cover. Then $X\cong S_1\times S_2/\langle (\gamma_1,\gamma_2)\rangle$,
with $\gamma_1$ the covering involution of the unramified double cover
$S_1\rightarrow T_1$ of an Enriques surface and $\gamma_2$ the
covering involution of some other (possibly ramified) double cover
$S_2\rightarrow T_2$.
\end{prp}

\begin{prf}
Let $\gamma$ be the covering involution of $\tilde{X}=S_1\times
S_2\rightarrow X$. Since
$\chi(\O_X)=\frac{1}{2}\chi(\O_{\tilde{X}})=2$ and $h^k(\O_X)\leq
h^k(\O_{\tilde{X}})$ we must have
$$h^k(\O_X)=\left\{ \begin{array}{ll}
  1 & \mbox{if }k=0\mbox{ or }4, \\
  0 & \mbox{otherwise.} \\
\end{array}\right.$$
In particular, the involution $\gamma$ does not fix either two-form
$p_i^*\sigma_i$, so it must act as $-1$.

Let $x_2\in S_2$. By integrating $p_1^*\sigma_1\bar{\sigma}_1$ over
$S_1\times x_2$ and $\gamma(S_1\times x_2)$ we see that $\gamma$ takes
$S_1\times x_2$ to $S_1\times y_2$ for some $y_2\in S_2$ (these are
calibrated submanifolds for $p_1^*\sigma_1\bar{\sigma}_1$). Similarly
$x_1\times S_2$ is sent to $y_1\times S_2$, and we conclude that
$\gamma=(\gamma_1,\gamma_2)$ acts independently on each
factor. If $\gamma_i$ both have fixed-points $p_i$, then
$(\gamma_1,\gamma_2)$ would have a fixed-point $(p_1,p_2)$, which is
impossible. Therefore $\gamma_1$, say, is fixed-point free and the
quotient $S_1/\langle\gamma_1\rangle$ is an Enriques surface.
\end{prf}

\subsection{An abelian fibration on $\tilde{X}$}

Although we want to understand fibrations $X\rightarrow\P^2$, it is
sometimes easier to work with the space $\tilde{X}$ because of the
product decompositions it admits in the different classes. Of course
$\tilde{X}$ is also
a fibration over $\P^2$ by virtue of the composition of maps
$\tilde{X}\rightarrow X\rightarrow\P^2$. Let $\Gamma$ be the order $d$
covering group of
$\tilde{X}\rightarrow X$; then $\Gamma$ is isomorphic to the quotient group 
$\pi_1(X)/\pi_1(\tilde{X})$. Let $t$ be a generic point of $\P^2$. The
generic fibre $A:=\pi^{-1}(t)$ of $\pi:X\rightarrow\P^2$ is an abelian
surface, and hence the generic fibre $\tilde{A}:=\tilde{\pi}^{-1}(t)$
of $\tilde{\pi}:\tilde{X}\rightarrow\P^2$ will be an unramified cover
of $A$ with covering group $\Gamma$. It follows that $\tilde{A}$ must
be a disjoint union of abelian surfaces. Let $\Gamma_0\subset\Gamma$
be the subgroup of elements which do {\em not\/} permute the connected
components of $\tilde{A}$. A connected component $\tilde{A}_0$ of
$\tilde{A}$ will be a cover of $A$ with covering group $\Gamma_0$;
thus $\Gamma_0$ is isomorphic to the quotient group
$\pi_1(A)/\pi_1(\tilde{A}_0)$, and in particular it is abelian. 

We now take the Stein factorization of
$\tilde{\pi}:\tilde{X}\rightarrow\P^2$, which gives
$$\begin{array}{ccc}
  \tilde{X} & \stackrel{\phi}{\rightarrow} & Z \\
  \downarrow & & \downarrow\tau \\
  X & \rightarrow & {\P}^2 \\
\end{array}$$
where the generic fibre of $\phi$ is an abelian surface, such as
$\tilde{A}_0$. Moreover, $\tau:Z\rightarrow\P^2$ is a Galois cover
with Galois group $\Gamma/\Gamma_0$; the points of $\tau^{-1}(t)$
correspond to the different connected components of $\tilde{A}$, and
$\Gamma/\Gamma_0$ acts freely and transitively on these components. In
particular, $Z\rightarrow\P^2$ is a $d^{\prime}$-to-1 cover, where
$d^{\prime}$ divides $d$. (This illustrates the utility of
the restrictions that $d$ is one in classes 5 and 6, and one or two
in class 3.)

The problem of constructing examples can be formulated as
follows. Given $\tilde{X}$ from one of the six classes, we look for a
fibration of $\tilde{X}$ by abelian surfaces over a surface $Z$. We
try to express $Z$ as a Galois cover $\tau:Z\rightarrow\P^2$, with
Galois group $\Gamma/\Gamma_0$. Of course $\tau$ must be ramified
since $\P^2$ is simply connected. There will be a quotient
$\tilde{X}/\Gamma_0$ which is also fibred by abelian surfaces
over $Z$. The problem is to extend the action of $\Gamma/\Gamma_0$ on
$Z$, which has fixed-points, to a fixed-point free action on
$\tilde{X}/\Gamma_0$; for then the quotient
$$(\tilde{X}/\Gamma_0)/(\Gamma/\Gamma_0)\cong \tilde{X}/\Gamma$$
is the required four-fold $X$, fibred by abelian surfaces over
$\P^2$. The following diagram summarizes this construction.
$$\begin{array}{ccc}
  \tilde{X} & & \\
  _{^{\Gamma_0}}\downarrow & \searrow & \\
  \tilde{X}/\Gamma_0 & \rightarrow & Z \\
  _{^{\Gamma/\Gamma_0}}\downarrow & & \downarrow _{^{\Gamma/\Gamma_0}} \\
  X & \rightarrow & {\P}^2 \\
\end{array}$$

If some element of $\Gamma/\Gamma_0$ fixes $s\in Z$, then that element
must act in a fixed-point free manner on the abelian surface fibre
above $s$. Usually the quotient of this action will again be an
abelian surface, and thus the element should act by a
translation. This suggests that $\Gamma/\Gamma_0$ should be
abelian, and we will keep this as a guiding principle when
constructing examples.

However, there may not be a point $s\in Z$ fixed by the entire group
$\Gamma/\Gamma_0$, so we won't find necessarily find that the entire
group acts as translations on a single abelian surface. Moreover, $s$
could lie in the discriminant locus of the abelian surface
fibration, implying that the fibre above $s$ is a degeneration rather
than a smooth abelian surface.


\subsection{Direct image sheaves}

Returning to $\pi:X\rightarrow\P^2$, we have the following results
concerning direct images of the structure sheaf.

\begin{lem}
The sheaves $R^0\pi_*\O_X$, $R^1\pi_*\O_X$, and $R^2\pi_*\O_X$ on 
$\P^2$ are locally free.
\end{lem}

\begin{prf}
Theorem 2.1 of Koll{\'a}r~\cite{kollar86i} says that 
$R^i\pi_*K_X$ is torsion-free for $i\geq 0$ (and vanishes for
$i>2$), but $K_X\cong\O_X$ since $X$ has trivial canonical
bundle. Corollary 3.9 of Koll{\'a}r~\cite{kollar86ii} then says that
in our case $R^i\pi_*\O_X$ is reflexive, but reflexive sheaves on
surfaces are locally free. Note that $R^1\pi_*\O_X$ has rank two,
while $R^0\pi_*\O_X$ and $R^2\pi_*\O_X$ are line bundles.
\end{prf}

\begin{lem}
We have $R^0\pi_*\O_X\cong\O_{\P^2}$ and
$R^2\pi_*\O_X\cong\O_{\P^2}(-3)$.
\end{lem}

\begin{prf}
Since $\pi:X\rightarrow\P^2$ has connected fibres,
$R^0\pi_*\O_X\cong\O_{\P^2}$ is automatic. Since $X$ has trivial
canonical bundle, we have $h^4(X,\O_X)=1$. This number can also be
computed using the Leray spectral sequence for $X\rightarrow\P^2$
$$E^{p,q}_2:=\H^q(\P^2,R^p\pi_*\O_X)\Rightarrow\H^{p+q}(X,\O_X).$$
Thus we have
$$\H^4(X,\O_X)\cong\H^2(\P^2,R^2\pi_*\O_X)\cong\H^0(\P^2,(R^2\pi_*\O_X)^*\otimes K_{\P^2})^*$$
implying that $R^2\pi_*\O_X\cong\O_{\P^2}(-3)$.
\end{prf}

\begin{lem}
Let $V$ denote the rank two bundle $R^1\pi_*\O_X$. Then
$$\H^q(\P^2,V)\cong \H^{q+1}(X,\O_X)$$
for all $q$.
\end{lem}

\begin{prf}
The previous lemma implies that the Leray spectral sequence
degenerates at $E^{p,q}_2$, and hence
$$\H^q(\P^2,R^p\pi_*\O_X)\cong \H^{p+q}(X,\O_X).$$
\end{prf}

\begin{rmk}
Recall that
$$h^i(X,\O_X)=h^0(X,\Omega^i_X)\leq
h^0(\tilde{X},\Omega^i_{\tilde{X}})$$
because any holomorphic form on $X$ must lift to $\tilde{X}$. The
right hand side can easily be calculated in each of the six
classes. We also have
$$\chi(\O_{\tilde{X}})=d\chi(\O_X)$$
and $h^3(X,\O_X)=h^1(X,\O_X)$ because of Serre duality and the
triviality of the canonical bundle.

There are therefore a limited number of possibilities for
$h^i(X,\O_X)$ in each of the six classes, and these lead to a limited
number of possibilities for the rank two bundle $V$, which
we will explore in the following sections.
\end{rmk}

\section{$\tilde{X}$ an abelian four-fold}

There are four possibilities for
$(h^q(\P^2,V))_{q=0}^2=(h^i(X,\O_X))_{i=1}^3$, namely
$$(4,6,4),\qquad (3,4,3),\qquad (2,2,2),\qquad\mbox{or}\qquad
(1,0,1).$$

\begin{lem}
\label{abelian_fourfold}
In the first subcase $X$ is an abelian four-fold. However, an abelian
four-fold $X$ cannot admit a fibration by abelian surfaces over $\P^2$.
\end{lem}

\begin{prf}
If $h^0(X,\Omega^1_X)=h^1(X,\O_X)=4$, then all holomorphic one-forms
on $\tilde{X}$ must be invariant under the action of the covering
group of $\tilde{X}\rightarrow X$, so that they descend to
$X$. Moreover, the Albanese variety $A$ of $X$ will be an abelian
four-fold and the composition
$$\tilde{X}\rightarrow X\rightarrow A$$
of the covering and Albanese maps will be an isogeny. This forces
$\tilde{X}\rightarrow X$ to also be an isogeny, and hence $X$ is an
abelian four-fold.

Suppose there is a fibration on $X$ whose generic fibre is an
abelian surface $B\subset X$. Then $B$ must be the translation of an
abelian subgroup $B_0\subset X$, and the deformations of $B$ in $X$
will be parametrized by the quotient group $X/B_0$. In other words,
the base of the fibration on $X$ must be an abelian surface $X/B_0$,
{\em not\/} the projective plane $\P^2$.
\end{prf}

\begin{lem}
\label{343}
The second subcase is impossible: there does not exist a rank two
vector bundle $V$ on $\P^2$ with
$$(h^q(\P^2,V))_{q=0}^2=(3,4,3).$$
\end{lem}

\begin{prf}
By assumption $V$ admits three independent sections. If these sections
don't generate the fibre of $V$ at a generic point, then they define a
subbundle; but then $V$ will be a direct sum of line bundles,
contradicting $h^1(\P^2,V)=4$. Therefore the evaluation map
$$\H^0(\P^2,V)\otimes\O_{\P^2}\rightarrow V$$
is surjective at a generic point; let the kernel be $\O_{\P^2}(k)$
with $k\leq 0$ and let the cokernel be the torsion sheaf $\F$. We
have a resolution of $\F$
$$0\rightarrow\O_{\P^2}(k)\rightarrow \O_{\P^2}^{\oplus 3}\rightarrow
V\rightarrow\F\rightarrow 0$$
which can be used to compute the cohomology of $\F$. Many of the terms
will vanish; for instance
$$h^1(\O_{\P^2}(k))=h^2(\O_{\P^2}(k))=h^1(\O_{\P^2}^{\oplus
  3})=h^2(\O_{\P^2}^{\oplus 3})=0,$$
whereas $h^2(V)=3$ by hypothesis. This will imply that $h^2(\F)=3$, a
contradiction since $\F$ is a torsion sheaf supported on a subscheme
  of $\P^2$ of dimension at most one.
\end{prf}

\begin{lem}
\label{222}
The third subcase is impossible: there does not exist a rank two
vector bundle $V$ on $\P^2$ with
$$(h^q(\P^2,V))_{q=0}^2=(2,2,2).$$
\end{lem}

\begin{prf}
By assumption $V$ admits two independent sections. As in the proof of
Lemma~\ref{343}, the evaluation map
$$\H^0(\P^2,V)\otimes\O_{\P^2}\rightarrow V$$
is surjective at a generic point; let the cokernel be the torsion
sheaf $\F$. Then we have a short exact sequence
$$0\rightarrow\O_{\P^2}^{\oplus 2}\rightarrow
V\rightarrow\F\rightarrow 0$$
whose corresponding long exact sequence yields
$$\H^2(\O_{\P^2}^{\oplus 2})\rightarrow \H^2(V)\rightarrow
\H^2(\F)\rightarrow 0.$$
The left-most term vanishes, implying that
$$h^2(\F)=h^2(V)=2.$$
This is a contradiction since $\F$ has support of dimension at
most one.
\end{prf}

\begin{lem}
\label{101}
In the fourth subcase $V\cong\O_{\P^2}\oplus\O_{\P^2}(-3)$.
\end{lem}

\begin{prf}
By assumption $V$ admits a section. Let the cokernel of
the evaluation map
$$\H^0(\P^2,V)\otimes\O_{\P^2}\rightarrow V$$
be $\O_{\P^2}(k)\otimes{\I}_Z$, where $k=c_1(V)$ and $\I_Z$ is the ideal
sheaf of a zero-dimensional subscheme $Z\subset\P^2$. Then we have a
short exact sequence
$$0\rightarrow\O_{\P^2}\rightarrow
V\rightarrow\O_{\P^2}(k)\otimes{\I}_Z\rightarrow 0$$
whose corresponding long exact sequence yields
$$\H^2(\O_{\P^2})\rightarrow \H^2(V)\rightarrow
\H^2(\O_{\P^2}(k)\otimes{\I}_Z)\rightarrow 0.$$
The left-most term vanishes, implying that
$$h^2(\O_{\P^2}(k))=h^2(\O_{\P^2}(k)\otimes\I_Z)=h^2(V)=1.$$
Therefore $c_1(V)=k=-3$.

Next the Riemann-Roch formula gives
$$\chi(V)=2\chi(\O_{\P^2})+\frac{c_1(V)(c_1(V)+3)}{2}-c_2(V)$$
and substituting $\chi(V)=2$ and $c_1(V)=-3$ yields $c_2(V)=0$. It
follows that the zero-dimensional subscheme $Z$ is empty. Moreover,
the short exact sequence
$$0\rightarrow\O_{\P^2}\rightarrow
V\rightarrow\O_{\P^2}(-3)\rightarrow 0$$
must split, completing the proof.
\end{prf}

We can construct an example in this fourth subcase.

\begin{exm}
Recall the discussion of Subsection~2.2. The abelian surface fibration
on $X$ induces an abelian surface fibration on $\tilde{X}$, with base
$Z$. The argument in the proof of Lemma~\ref{abelian_fourfold} shows
that $Z$ must be an abelian surface. First we need to find a Galois
cover $Z\rightarrow\P^2$ with $Z$ an abelian surface. Then we need to
lift the group action to a fixed-point free action on $\tilde{X}$.

Let $E$ be an elliptic curve and let $Z$ be $E\times E$. Since $E$ is
a double cover of $\P^1$, we get the composition of maps
$$Z=E\times E\rightarrow\P^1\times\P^1\rightarrow\mathrm{Sym}^2\P^1\cong\P^2.$$
This expresses $Z$ as a Galois cover of
$\P^2$. Let $z_1$ and $z_2$ be coordinates on the first and second
copies of $E$ respectively. Then the Galois group $\Gamma$ is
generated by elements $\gamma_1$, $\gamma_2$, and $\gamma_3$ which act
by
$$\gamma_1.(z_1,z_2)=(-z_1,z_2),$$
$$\gamma_2.(z_1,z_2)=(z_1,-z_2),$$
$$\mbox{and}\qquad\gamma_3.(z_1,z_2)=(z_2,z_1).\qquad\phantom{and}$$
Note that $\gamma_3\gamma_1\gamma_3=\gamma_2$ and
$\Gamma$ is isomorphic to the dihedral group $D_8$: for instance, the
element $\gamma_1\gamma_3$ corresponds to a rotation of order four.

Next we lift the $\Gamma$-action to the abelian four-fold
$\tilde{X}$. For simplicity, we assume that $\tilde{X}$ is $E\times
E\times E_3\times E_4$ where $E_3$ and $E_4$ are arbitrary elliptic
curves. Let $z_3$ and $z_4$ be coordinates on $E_3$ and $E_4$
respectively, and write $E_3\cong\C/\langle 1,\tau_3\rangle$ and
$E_4\cong\C/\langle 1,\tau_4\rangle$. We define a $\Gamma$-action on
$E_3\times E_4$ by
$$\gamma_1.(z_3,z_4)=(z_3+\frac{1}{2},-z_4+\frac{1}{4}),$$
$$\gamma_2.(z_3,z_4)=(z_3+\frac{1}{2},-z_4+\frac{3}{4}),$$
$$\mbox{and}\qquad\gamma_3.(z_3,z_4)=(z_3+\frac{1}{2}\tau_3,-z_4).\qquad\phantom{and}$$
The
`rotation' $\gamma_1\gamma_3$ acts by a translation
$$\gamma_1\gamma_3.(z_3,z_4)=(z_3+\frac{1}{2}+\frac{1}{2}\tau_3,z_4+\frac{1}{4})$$
of order four. Thus the quotient by $\langle\gamma_1\gamma_3\rangle$
is an abelian surface, and further quotienting by $\Z_2\cong
\Gamma/\langle\gamma_1\gamma_3\rangle$ produces a bielliptic
surface. Note that every non-trivial element of $\Gamma$ acts in a
fixed-point free manner on $E_3\times E_4$, and hence the product
$\Gamma$-action on $\tilde{X}=E\times E\times E_3\times E_4$ is also
fixed-point free. This action preserves the four-form $dz_1\wedge
dz_2\wedge dz_3\wedge dz_4$; the only other forms it preserves are the
one-form $dz_3$ and the three-form $dz_1\wedge dz_2\wedge dz_4$. Thus
the quotient $X:=\tilde{X}/\Gamma$ has trivial canonical bundle and
$$(h^i(X,\O_X))_{i=1}^3=(1,0,1).$$

If $t$ is a generic point of $\P^2$, then the fibre of
$Z\rightarrow\P^2$ over $t$ consists of eight points, and hence the
fibre of $\tilde{X}\rightarrow\P^2$ over $t$ consists of eight copies
of $E_3\times E_4$. Since $\Gamma$ acts freely and transitively on
these copies of $E_3\times E_4$, the fibre of $X\rightarrow\P^2$ over
$t$ will consist of a single copy of $E_3\times E_4$. Thus $X$ is a
generically isotrivial abelian surface fibration over $\P^2$: generic
fibres are all isomorphic to the fixed abelian surface $E_3\times
E_4$.
\end{exm}

\section{$\tilde{X}=A\times S$, $A$ and $S$ abelian and K3 surfaces}

There are two possibilities for
$(h^q(\P^2,V))_{q=0}^2=(h^i(X,\O_X))_{i=1}^3$, namely
$$(2,2,2)\qquad\mbox{or}\qquad
(1,0,1).$$
By Lemma~\ref{222} there is no rank two vector bundle $V$ on $\P^2$ in
the first subcase. In particular, if $X$ is $A\times S$ itself, then
$X$ cannot be fibred by abelian surfaces over $\P^2$. In the second
subcase $V\cong\O_{\P^2}\oplus\O_{\P^2}(-3)$ by Lemma~\ref{101} and
$X$ must be a non-trivial quotient of $\tilde{X}=A\times S$.

We can construct an example in this second subcase.

\begin{exm}
Recall again the discussion in Subsection~2.2. There should be an
induced abelian surface fibration on $\tilde{X}=A\times S$. An obvious
choice would be the trivial fibration $\tilde{X}\rightarrow S$. Thus
we need to express $S$ as a Galois cover of $\P^2$, so we let $S$ be a
hyperelliptic K3 surface, i.e., a double cover of $\P^2$ branched over
a sextic curve. If $\sigma$ is the holomorphic symplectic two-form on
$S$, then $\gamma^*\sigma=-\sigma$ where $\gamma$ is the covering
involution of $S\rightarrow\P^2$.

Next we extend the action of $\gamma$ to an involution on
$\tilde{X}=A\times S$. The action on $A$ must take the two-form
$dz_1\wedge dz_2$ to $-dz_1\wedge dz_2$, so that the involution
preserves the triviality of the canonical bundle on $\tilde{X}$. Thus
we let $A$ be the product of two elliptic curves $E_1\cong\C/\langle
1, \tau_1\rangle$ and $E_2\cong\C/\langle 1, \tau_2\rangle$ with
coordinates $z_1$ and $z_2$ respectively, and we let $\gamma$ act by
$$\gamma.(z_1,z_2)=(z_1+\frac{1}{2},-z_2).$$
The quotient $E_1\times E_2/\langle\gamma\rangle$ is a bielliptic
surface. This involution on $E_1\times E_2$ is fixed-point
free, and therefore the product action of $\gamma$ on $\tilde{X}=E_1\times
E_2\times S$ is also fixed-point free. The action preserves the
four-form $dz_1\wedge dz_2\wedge \sigma$, the one-form $dz_1$ and the
three-form $dz_2\wedge\sigma$. Thus the quotient
$X:=\tilde{X}/\langle\gamma\rangle$ has trivial canonical bundle and
$$(h^i(X,\O_X))_{i=1}^3=(1,0,1).$$

Over a generic point $t\in\P^2$, the fibre of $S\rightarrow\P^2$
consists of two points, the fibre of $\tilde{X}\rightarrow\P^2$
consists of two copies of $E_1\times E_2$, and hence the fibre of
$X\rightarrow\P^2$ will consist of a single copy of $E_1\times
E_2$. This expresses $X$ as a generically isotrivial abelian surface
fibration over $\P^2$.
\end{exm}

\section{$\tilde{X}=S_1\times S_2$, $S_1$ and $S_2$ K3 surfaces}

There are two possibilities for
$(h^q(\P^2,V))_{q=0}^2=(h^i(X,\O_X))_{i=1}^3$, namely
$$(0,2,0)\qquad\mbox{or}\qquad
(0,0,0).$$
In the first subcase $d$ must equal one, so that
$X=\tilde{X}=S_1\times S_2$. In the second subcase $d$ must equal two,
so that $\tilde{X}$ is a double cover of $X$.

The next two lemmas could have been proved earlier, but we did not
need them until now; first some notation. Let
$\Delta\subset\P^2$ be the degeneracy locus which parametrizes
singular fibres of $\pi:X\rightarrow\P^2$, and let
$U=\P^2\backslash\Delta$ be its complement. For $t\in U$ we have a
canonical isomorphism  
$$\Lambda^2(R^1\pi_*\O_X)_t=\Lambda^2\H^1(X_t,\O_{X_t})\cong\H^2(X_t,\O_{X_t})=(R^2\pi_*\O_X)_t$$
and thus
$$\Lambda^2(R^1\pi_*\O_X)|_U\cong R^2\pi_*\O_X|_U.$$

\begin{lem}
\label{first_chern}
There is an inclusion
$$\Lambda^2(R^1\pi_*\O_X)\hookrightarrow R^2\pi_*\O_X$$
which extends the above isomorphism. In particular, we have
$$c_1(V)=c_1(R^1\pi_*\O_X)\leq c_1(R^2\pi_*\O_X)=-3.$$
Moreover, if either
\begin{itemize}
\item $\H^1(X_t,\C)\geq 2$ where $X_t$ is a generic singular fibre of
  $X\rightarrow\P^2$
\item or the monodromy of $R^1\pi_*\O_X$ around $\Delta$ is unipotent
\end{itemize}
then the inclusion above is actually an isomorphism and hence
$c_1(V)=-3$.
\end{lem}

\begin{prf}
The proofs of these statements can be found in
Matsushita~\cite{matsushita05}, particularly Proposition~2.13 and
Section 3. For instance, the inclusion
$$\Lambda^2(R^1\pi_*\O_X)\hookrightarrow R^2\pi_*\O_X$$
is a consequence of functoriality of the variation of Hodge
structure.
\end{prf}

\begin{lem}
\label{kollar}
The higher cohomology of $V(k):=V\otimes\O_{\P^2}(k)$ vanishes for $k$
positive, i.e.,
$$h^q(\P^2,V(k))=0$$
for $k\geq 1$ and $q\geq 1$.
\end{lem}

\begin{prf}
This is the third part of Theorem~2.1 in Koll{\'a}r~\cite{kollar86i}.
\end{prf}

We now return to the two subcases of this section.

\begin{lem}
\label{020}
The first subcase is impossible: there does not exist a rank two
vector bundle $V$ on $\P^2$ with
$$(h^q(\P^2,V))_{q=0}^2=(0,2,0),$$
$c_1(V)\leq -3$, and $h^q(\P^2,V(1))=0$ for $q\geq 1$.
\end{lem}

\begin{prf}
The Riemann-Roch formula gives
$$\chi(V)=2\chi(\O_{\P^2})+\frac{c_1(V)(c_1(V)+3)}{2}-c_2(V)$$
whereas we know that $\chi(V)=-2$. Therefore
\begin{eqnarray*}
\chi(V(1)) & = &
2\chi(\O_{\P^2})+\frac{(c_1(V)+2)(c_1(V)+2+3)}{2}-(c_2(V)+c_1(V)+1) \\
 & = & \left(
2\chi(\O_{\P^2})+\frac{c_1(V)(c_1(V)+3)}{2}-c_2(V)\right)+c_1(V)+4 \\
 & = & \chi(V)+c_1(V)+4 \\
 & = & c_1(V)+2.
\end{eqnarray*}
Then $c_1(V)\leq -3$ (from Lemma~\ref{first_chern}) implies
$$\chi(V(1))\leq -1$$
whereas $h^q(\P^2,V(1))=0$ for $q\geq 1$ (from Lemma~\ref{kollar})
implies the contradictory
$$\chi(V(1))=h^0(\P^2,V(1))\geq 0.$$
\end{prf}

\begin{cor}
\label{impossible1}
The four-fold $S_1\times S_2$ given by the product of two K3 surfaces
cannot be fibred by abelian surfaces over $\P^2$.
\end{cor}

In the second subcase, recall that $X=S_1\times
S_2/\langle(\gamma_1,\gamma_2)\rangle$ by Proposition~\ref{Enriques},
where $\gamma_1$ is the covering involution of the double cover
$S_1\rightarrow T_1$ of an Enriques surface and $\gamma_2$ is the
covering involution of some other (possibly ramified) double cover
$S_2\rightarrow T_2$.

\begin{prp}
\label{impossible2}
The four-fold $X=S_1\times S_2/\langle(\gamma_1,\gamma_2)\rangle$
cannot be fibred by abelian surfaces over $\P^2$.
\end{prp}

\begin{prf}
Suppose we have a fibration $X\rightarrow B$ by abelian surfaces, with
$B$ some surface, and let $A$ be a generic fibre. By composing the
inclusion $A\hookrightarrow X$, the (possibly ramified) double cover
$$X=S_1\times S_2/\langle(\gamma_1,\gamma_2)\rangle\rightarrow
T_1\times T_2$$
and the projection to $T_1$ we get a map $A\rightarrow T_1$. The
Enriques surface has Picard number $\rho(T_1)=10$, whereas the abelian
surface has Picard number $\rho(A)$ at most four. Therefore
$A\rightarrow T_1$ cannot be surjective. 

If the image of $A\rightarrow T_1$ were zero-dimensional it would have
to be a single point as $A$ is connected. Then the deformations of $A$
in $X$ would be parametrized either by $T_1$ or by a double cover of
$T_1$. In other words, $B$ would be isomorphic either to $T_1$ or to a
double cover of $T_1$; but neither of these is $\P^2$ (in the latter
case, the Picard number $\rho(B)\geq \rho(T_1)=10$).

Finally, suppose the image of $A\rightarrow T_1$ is a curve $C\subset
T_1$. The base of a fibration $A\rightarrow C$ of an abelian surface
over a curve must have genus zero or one. If $C$ were a rational curve
it would be rigid in $T_1$; but as we move $A$ in $X$ the curve $C$
should also move in $T_1$. Therefore $C$ is an elliptic curve, a
fibre of an elliptic fibration $T_1\rightarrow\P^1$ on the Enriques
surface. This induces a map $B\rightarrow\P^1$ (the point in $B$
corresponding to the fibre $A$ is mapped to the point in $\P^1$
corresponding to the fibre $C$), and therefore the Picard number
$\rho(B)\geq 2$. Again we see that $B$ cannot be $\P^2$.
\end{prf}

\section{$\tilde{X}=E\times Y$, $E$ an elliptic curve, $Y$ a
  Calabi-Yau three-fold}

There is only one possibility for
$(h^q(\P^2,V))_{q=0}^2=(h^i(X,\O_X))_{i=1}^3$, namely
$(1,0,1)$. However, we can't conclude anything about the degree $d$ of
the cover $\tilde{X}\rightarrow X$, so $X$ itself does not necessarily
decompose as $E\times Y$. By Lemma~\ref{101} the bundle $V$ must be
isomorphic to $\O_{\P^2}\oplus\O_{\P^2}(-3)$.

\begin{exm}
An obvious example is given by choosing the Calabi-Yau three-fold $Y$
to be an elliptic fibration over $\P^2$, and letting
$X=\tilde{X}=E\times Y$. Then $X$ is an abelian surface fibration over
$\P^2$ whose fibres are products of $E$ and the elliptic fibres of
$\eta:Y\rightarrow\P^2$. Let us examine $Y$ more closely.

By assumption $X\rightarrow\P^2$ has a section $s$, and therefore
$Y\rightarrow\P^2$ must admit a section too, which we also denote by
$s$. Some singular fibres of $Y\rightarrow\P^2$ may consist of several
irreducible components; by blowing down all components of the singular
fibres that the section does not meet, we obtain the Weierstrass model
$\bar{\eta}:\overline{Y}\rightarrow\P^2$ of $Y$. We recall the
description of the Weierstrass model $\overline{Y}$ as a family of
plane cubics over $\P^2$ (for example, Friedman~\cite{friedman98}
describes the theory of elliptic surfaces, but his arguments on pages
180-181 apply to general elliptic fibrations).

A single elliptic curve $E$ is embedded as a cubic in the plane by
$$E\hookrightarrow\P(\mathrm{H}^0(E,\O(3p))^*)$$
where $p$ is some basepoint on $E$. In the relative version the
basepoint in each curve is given by the section $s$, and the embedding
looks like
$$\overline{Y}\hookrightarrow\P(W^*)$$
where $W$ is the rank three bundle
$\bar{\eta}_*\O_{\overline{Y}}(3s)$. One can show that
$$W\cong L^{-2}\oplus L^{-3}\oplus\O_{\P^2}$$
where $L$ is the dual of the line bundle
$R^1\bar{\eta}_*\O_{\overline{Y}}$. Since $Y$ (and hence
$\overline{Y}$) has trivial canonical bundle, we must have
$$R^1\bar{\eta}_*\O_{\overline{Y}}=R^1\eta_*\O_Y\cong K_{\P^2}=\O_{\P^2}(-3)$$
and therefore $L\cong\O_{\P^2}(3)$ and
$$W^*\cong\O_{\P^2}(6)\oplus\O_{\P^2}(9)\oplus\O_{\P^2}.$$
Let $(x,y,z)$ denote fibre coordinates on $W^*$ corresponding to the above 
decomposition. Then $\overline{Y}$ is a family of
cubics in $\P(W^*)$ over the projective plane given by
the equation
$$y^2z=x^3+axz^2+bz^3,$$
where $a$ and $b$ are sections of $\O_{\P^2}(12)$ and $\O_{\P^2}(18)$
respectively. The discriminant locus $\Delta\subset\P^2$ parametrizing
singular fibres is the zero locus of
$$4a^3+27b^2\in\H^0(\P^2,\O(36)).$$
For generic $a$ and $b$, $\Delta$ is a smooth curve in the plane and
every singular fibre is a nodal rational curve (in particular, this
means that the Weierstrass model $\overline{Y}$ is smooth and
isomorphic to $Y$ itself). It looks like the family of such fibrations
$\overline{Y}\rightarrow\P^2$ will have dimension
$$h^0(\P^2,\O(12))+h^0(\P^2,\O(18))=13+19=32,$$
though actually $a$ and $b$ are only defined up to a rescaling
$$(a,b)\mapsto (\lambda^4a,\lambda^6b)$$
where $\lambda\in\mathbb{C}^*$, and we also need to quotient by the
automorphism group $\mathrm{PGL}(3,\mathbb{C})$ of $\P^2$. Thus we are
left with a family of dimension
$$32-1-8=23.$$
In addition, the moduli space of elliptic curves is one-dimensional,
so $X=E\times Y$ depends on $24$ parameters.
\end{exm}

\begin{rmk}
Suppose that $X=E\times Y$. The image of $Y_t:=\{t\}\times Y$ under
the projection $\pi:X\rightarrow\P^2$ must be at least
one-dimensional. If the image is a curve $C\subset\P^2$, then the
(trivial) normal bundle of $Y_t$ in $E\times Y$ would be the pull-back
of the non-trivial normal bundle of $C$ in $\P^2$, a
contradiction. Therefore $Y_t\rightarrow\P^2$ is surjective. The
fibres of $Y_t\rightarrow\P^2$ must be connected, since the fibres of
$X\rightarrow\P^2$ are connected; this makes $Y_t$ into an elliptic
fibration over $\P^2$. Moreover, the morphism to $\P^2$ is canonically
given by
$$Y_t\rightarrow\P(\H^0(Y_t,L_t)^*)$$
where $L_t$ is the pull-back of the hyperplane line bundle from
$\P^2$. Since $Y_t\cong Y$ for all $t\in E$, and the line bundle $L_t$
cannot be deformed because $\H^1(Y,\O)=0$, the elliptic fibration
$Y_t\rightarrow\P^2$ will not depend on $t$ either. Therefore
$X=E\times Y\rightarrow\P^2$ is exactly as constructed in the example,
the product of an elliptic curve and an elliptic fibration over
$\P^2$.

The general situation, when $X$ is a quotient of $\tilde{X}=E\times
Y$, is more complicated.
\end{rmk}

%

\section{$\tilde{X}=X$ a Calabi-Yau four-fold}

There is one possibility for
$(h^q(\P^2,V))_{q=0}^2=(h^i(X,\O_X))_{i=1}^3$, namely
$(0,0,0)$. By Lemma~\ref{first_chern} we know that $c_1(V)\leq
-3$, and by Lemma~\ref{kollar} we know that $h^q(\P^2,V(k))=0$ for
$k\geq 1$ and $q\geq 1$.

\begin{lem}
\label{000}
The bundle $V$ on $\P^2$ is isomorphic to either
$\O_{\P^2}(-1)\oplus\O_{\P^2}(-2)$ or
$\O_{\P^2}(-2)\oplus\O_{\P^2}(-2)$. Since the latter bundle has first
Chern class $-4$, it cannot occur if either of the hypotheses 
of Lemma~\ref{first_chern} regarding the singular fibres of
$X\rightarrow\P^2$ are satisfied (i.e., first cohomology at least
two-dimensional or unipotent monodromy).
\end{lem}

\begin{prf}
As in the proof of Lemma~\ref{020}, we can use the Riemann-Roch
formula and $\chi(V)=0$ to calculate
\begin{eqnarray*}
\chi(V(1)) & = & \chi(V)+c_1(V)+4 \\
 & = & c_1(V)+4.
\end{eqnarray*}
Lemma~\ref{kollar} gives
$$\chi(V(1))=h^0(\P^2,V(1))\geq 0$$
and hence $c_1(V)\geq -4$. This complements the inequality
$c_1(V)\leq -3$ of Lemma~\ref{first_chern}.

First suppose that $c_1(V)=-3$; then $h^0(\P^2,V(1))=1$ and we have a
morphism $\O_{\P^2}(-1)\rightarrow V$. The cokernel of this map is of
the form $\O_{\P^2}(k)\otimes\I_Z$, where $\I_Z$ is the ideal sheaf of
some zero-dimensional subscheme $Z\subset\P^2$. Moreover $k=-2$
because $c_1(V)=-3$. Substituting $\chi(V)=0$ and $c_1(V)=-3$ into the
Riemann-Roch formula gives $c_2(V)=2$, which implies that the
zero-dimensional subscheme $Z$ is empty. The short exact sequence
$$0\rightarrow\O_{\P^2}(-1)\rightarrow
V\rightarrow\O_{\P^2}(-2)\rightarrow 0$$
must split, giving $V\cong\O_{\P^2}(-1)\oplus\O_{\P^2}(-2)$.

A similar argument works when $c_1(V)=-4$. Although $h^0(\P^2,V(1))$
is now zero, we can instead compute $h^0(\P^2,V(2))$ and find that it
equals two. Let the cokernel of the evaluation map
$$\H^0(\P^2,V(2))\otimes\O_{\P^2}(-2)\rightarrow V$$
be the torsion sheaf $\F$. The Riemann-Roch formula yields $c_2(V)=4$,
and since the Chern classes of $V$ agree with those of
$\O_{\P^2}(-2)\oplus\O_{\P^2}(-2)$ the evaluation map must be an
isomorphism and $\F$ must vanish. Therefore
$V\cong\O_{\P^2}(-2)\oplus\O_{\P^2}(-2)$. 
\end{prf}

We will describe here an example whose fibres are products of elliptic
curves; in Section~9 we will describe an example whose fibres are
Jacobians of genus two curves. Both examples will have
$V\cong\O_{\P^2}(-1)\oplus\O_{\P^2}(-2)$.

\begin{exm}
For $i=1$ and $2$, let $\pi_i:X_i\rightarrow\P^2$ be elliptic
fibrations with $R^1\pi_{i*}\O_{X_i}\cong\O_{\P^2}(-i)$. For example,
we could let $X_1$ be a Weierstrass model given by a family of cubics
$$y^2z=x^3+axz^2+bz^3$$
in the $\P^2$-bundle
$$\P(\O_{\P^2}(2)\oplus\O_{\P^2}(3)\oplus\O_{\P^2})$$
over the projective plane, where $a$ and $b$ are generic sections of
$\O_{\P^2}(4)$ and $\O_{\P^2}(6)$ respectively. By genericity, the
discriminant locus
$$\Delta_1=\{4a^3+27b^2=0\}\subset\P^2$$
will be a smooth curve; thus the singular fibres will all be nodal
rational curves and $X_1$ will be smooth. Likewise, we let $X_2$ be a
Weierstrass model given by
$$y^2z=x^3+cxz^2+dz^3$$
in
$$\P(\O_{\P^2}(4)\oplus\O_{\P^2}(6)\oplus\O_{\P^2}),$$
where $c$ and $d$ are generic sections of $\O_{\P^2}(8)$ and
$\O_{\P^2}(12)$ respectively. Again, the discriminant locus
$$\Delta_2=\{4c^3+27d^2=0\}\subset\P^2$$
will be a smooth curve, the singular fibres will all be nodal rational
curves, and $X_2$ will be smooth.

Next we define $X$ to be the fibre product $X_1\times_{\P^2}X_2$. Then
$X$ is a fibration over $\P^2$ whose fibres are products of elliptic
curves. The singular fibres of $X\rightarrow\P^2$ occur over
$\Delta_1\cup\Delta_2$; over $\Delta_1\backslash\Delta_2$ and
$\Delta_2\backslash\Delta_1$ these are products of nodal rational
curves and smooth elliptic curves, whereas over $\Delta_1\cap\Delta_2$
these are products of two nodal rational curves. By genericity, we can
assume that $\Delta_1$ and $\Delta_2$ intersect transversally (in a
collection of points $\Delta_1\cap\Delta_2$), and then a local
calculation shows that $X$ is smooth.

Now $\pi:X\rightarrow\P^2$ has
$$V=R^1\pi_*\O_X\cong R^1\pi_{1*}\O_{X_1}\oplus
R^1\pi_{2*}\O_{X_2}\cong\O_{\P^2}(-1)\oplus\O_{\P^2}(-2).$$
The Leray spectral sequence therefore implies that
$$h^k(\O_X)=\left\{ \begin{array}{ll}
  1 & \mbox{if }k=0\mbox{ or }4, \\
  0 & \mbox{otherwise.} \\
\end{array}\right.$$
We claim that the canonical bundle $K_X$ is trivial. Each fibre of
$X\rightarrow\P^2$ has trivial canonical bundle and trivial normal
bundle. Therefore the adjunction formula implies that $K_X$ is trivial
when restricted to any fibre; but then $K_X$ must be the pull-back
$\pi^*\O_{\P^2}(m)$ of a line bundle on $\P^2$. Since
$$h^0(K_X)=h^0(\Omega^4_X)=h^4(\O_X)=1,$$
the canonical bundle $K_X$ has a unique section up to scale. This
implies that $m=0$, proving the claim.

Finally, $X$ must be a Calabi-Yau four-fold as the only other
four-fold $X$ with trivial canonical bundle and $h^k(\O_X)$ as above
is the quotient of a product $S_1\times S_2$ of K3 surfaces by an
involution, but Proposition~\ref{impossible2} says that such a four-fold
cannot be fibred by abelian surfaces over $\P^2$.

Note that $a$, $b$, $c$, and $d$ are defined up to rescalings
$$(a,b)\mapsto (\lambda^4a,\lambda^6b)\qquad\mbox{and}\qquad
(c,d)\mapsto (\mu^4c,\mu^6d)$$
where $\lambda$ and $\mu$ are in $\C^*$. We also need to quotient by
the automorphism group of $\P^2$, and thus we get a family of examples
depending on
$$h^0(\P^2,\O(4))+h^0(\P^2,\O(6))+h^0(\P^2,\O(8))+h^0(\P^2,\O(12))-2-\mathrm{dim}\mathrm{PGL}(3,\C)\hspace*{5mm}$$
$$\hspace*{60mm}=5+7+9+13-2-8=24$$
parameters.
\end{exm}

\section{$\tilde{X}=X$ an irreducible holomorphic symplectic
  four-fold}

There is one possibility for
$(h^q(\P^2,V))_{q=0}^2=(h^i(X,\O_X))_{i=1}^3$, namely
$(0,1,0)$. Using arguments similar to those above, one can show that
$V$ is isomorphic to the cotangent bundle $\Omega^1$ on $\P^2$ (one
obtains a variation of the Euler sequence defining $\Omega^1$). In
fact, $\pi:X\rightarrow\P^2$ is a holomorphic Lagrangian fibration,
and Matsushita~\cite{matsushita05} showed more generally that if
$\pi:X\rightarrow\P^n$ is a Lagrangian fibration on an irreducible
holomorphic symplectic $2n$-fold $X$, then
$$R^i\pi_*\O_X\cong\Omega^i_{\P^n}$$
for all $i=0,\ldots,n$. The isomorphism $V\cong\Omega^1_{\P^2}$ is a
special case of this result.

Up to deformation, there are two known examples of irreducible
holomorphic symplectic four-folds. One is the Hilbert scheme
$\mathrm{Hilb}^2S$ of two points on a K3 surface $S$, which was first
described by Fujiki~\cite{fujiki83} as
$$\mathrm{Hilb}^2S=(\mathrm{Blow}_{\Delta}S\times S)/\Z_2.$$
The other is the generalized Kummer four-fold, due to
Beauville~\cite{beauville83}.

\begin{exm}
If $S$ is an elliptic K3 surface, then the elliptic fibration
$S\rightarrow\P^1$ induces a Lagrangian fibration
$$\mathrm{Hilb}^2S\rightarrow\mathrm{Sym}^2\P^1\cong\P^2$$
whose fibres are products of elliptic curves. Another example is given
by choosing $S$ to contain a genus two curve $C$, which will move in a
two-dimensional linear system $|C|\cong\P^2$. Denote this family of
curves by $\mathcal{C}\rightarrow\P^2$. Then the compactified relative
Jacobian
$$X:=\overline{\mathrm{Jac}}^0(\mathcal{C}/\P^2)\rightarrow\P^2$$
is both birational to and a deformation of $\mathrm{Hilb}^2S$. This
Lagrangian fibration by Jacobians of genus two curves is known as the
Beauville-Mukai integrable system~\cite{beauville99}.

Suppose that $S$ is a generic polarized K3 surface, i.e.,
$\mathrm{Pic}S$ is generated by the polarization $H$. Proposition~7.1
of Hassett and Tschinkel~\cite{ht99} states that the Hilbert scheme
$\mathrm{Hilb}^2S$ is a Lagrangian fibration by abelian surfaces over
$\P^2$ if and only if the polarization has degree $H^2=2m^2$ for some
integer $m>1$ (in the case $m=1$ and $H^2=2$, $\mathrm{Hilb}^2S$ is
obtained from the Beauville-Mukai integrable system via a Mukai
flop). Independently, Markushevich~\cite{markushevich06} and the
author~\cite{sawon07} gave an alternate proof of this result, and
generalized it to higher dimensions. In general, if $S$ is not
generic, $H^2=2m^2$ will only imply that $\mathrm{Hilb}^2S$ is
birational to a Lagrangian fibration.
\end{exm}

\begin{exm}
Debarre~\cite{debarre99} described Lagrangian fibrations on
generalized Kummer varieties (see also Example~3.8 in the author's
article~\cite{sawon03}). In particular, certain generalized Kummer
four-folds admit Lagrangian fibrations over $\P^2$ whose fibres are
abelian surfaces of polarization type $(1,3)$. Note that by
Proposition~5.3 of~\cite{sawon03} a generalized Kummer variety cannot
be fibred by principally polarized abelian varieties, such as
Jacobians of curves. Gulbrandsen~\cite{gulbrandsen07} and
Yoshioka~\cite{yoshioka06} investigated exactly which generalized
Kummer varieties are Lagrangian fibrations.
\end{exm}

\section{Fibrations by Jacobians}

In this final section we construct examples of four-folds with trivial
canonical bundles which are fibred over $\P^2$ by Jacobians of genus
two curves. Moreover, we classify all such four-folds. The results in
this section owe a lot to Markushevich's papers~\cite{markushevich95,
  markushevich96}; indeed, the examples we describe are the
Beauville-Mukai system (a holomorphic symplectic four-fold) and
Example~1 from~\cite{markushevich95} (which we prove is a Calabi-Yau
four-fold). Markushevich~\cite{markushevich96} classified
all irreducible holomorphic symplectic four-folds fibred by Jacobians;
our classification result is an extension of his, and uses similar
arguments.

We start by defining exactly what we mean by a fibration by
Jacobians. The following is essentially Definition~4
from~\cite{markushevich95}.

\begin{dfn}
Let $\mathcal{C}\rightarrow\P^2$ be a flat family of irreducible and
reduced genus two curves over the projective plane. We say that
$\mathcal{C}/\P^2$ has mild degenerations if
\begin{enumerate}
\item the total space of $\mathcal{C}$ is smooth,
\item the singular curves of the family have only nodes or cusps as
  singularities,
\item if $C_t$ is a singular curve with two distinct singular points
  $p_1$ and $p_2$, then the analytic germs $(\Delta_1,t)$ and
  $(\Delta_2,t)$ of the discriminant curves of the unfoldings of
  $(C_t,p_1)$ and $(C_t,p_2)$ induced by the family $\mathcal{C}/\P^2$
  must have different reduced tangent cones at $t$. 
\end{enumerate}
\end{dfn}

\begin{rmk}
The smoothness of the total space of $\mathcal{C}$ implies that a
generic singular curve will have a single node, and will occur in
codimension one in the family. In codimension two we may get curves
with cusps or a pair of nodes. The third condition means the
following. If $C_t$ has nodes at $p_1$ and $p_2$ then there will be
two branches of the discriminant curve $\Delta\subset\P^2$ at
$t$. Along one branch, the curves will have nodes whose limit is $p_1$
as we approach $C_t$, and along the other branch the curves will have
nodes whose limit is $p_2$. The mild degenerations condition requires
these two branches to meet transversally at $t$.
\end{rmk}

\begin{dfn}
By a fibration by Jacobians over $\P^2$ we will mean the compactified
relative Jacobian
$$X:=\overline{\mathrm{Jac}}^0(\mathcal{C}/\P^2)$$
of a family $\mathcal{C}/\P^2$ of irreducible and reduced genus two
curves with mild degenerations.
\end{dfn}

\begin{rmk}
Theorem~2 of~\cite{markushevich95} states that if $\mathcal{C}/\P^2$
has mild degenerations then its compactified relative Jacobian $X$ is
smooth and its canonical bundle $K_X$ is trivial when restricted to
the fibres of $X\rightarrow\P^2$. Markushevich strengthened this
theorem in~\cite{markushevich96} by giving necessary and sufficient
conditions on a family of curves for its compactified relative
Jacobian to be smooth. We will stick with the mild degenerations
condition for simplicity.
\end{rmk}

\begin{lem}
Let $\pi:X\rightarrow\P^2$ be a fibrations by Jacobians over $\P^2$
constructed from a family of curves $\mathcal{C}/\P^2$. Let $Y$ denote
the total space of $\mathcal{C}/\P^2$. Then $Y$ is a branched double
cover of the $\P^1$-bundle $\P(R^1\pi_*\O_X)$ over $\P^2$.
\end{lem}

\begin{prf}
A genus two curve $C$ is hyperelliptic: the canonical map
$$C\rightarrow\P(\H^0(C,K_C)^*)\cong\P(\H^1(C,\O_C))$$
expresses $C$ as a double cover of $\P^1$, branched over six
points. The relative version of this expresses $Y$ as a double cover
of the $\P^1$-bundle $\P(R^1f_*\O_Y)$ over $\P^2$, branched over a
divisor which intersects each fibre $\P^1$ in six points, counted with
multiplicity. Here $f$ denotes the projection $Y\rightarrow\P^2$.

It remains to show that
$$R^1f_*\O_Y\cong R^1\pi_*\O_X,$$
but this is Lemma~5 in~\cite{sawon08iv}. More precisely, Lemma~5 is
for Lagrangian fibrations, but the same argument works provided that
the generic singular curves of $\mathcal{C}/\P^2$ have only single
nodes, which is guaranteed by the mild degenerations hypothesis.
\end{prf}

We now reverse the construction by starting with a rank two bundle $W$
on $\P^2$. Let $Y$ be a double cover of the $\P^1$-bundle
$h:\P(W)\rightarrow\P^2$ branched over a divisor corresponding to a
generic section of
$$\O_{\P(W)}(6)\otimes h^*\O(2d).$$
Such a divisor will intersect each fibre $\P^1$ in six points, counted
with multiplicity, and hence $f:Y\rightarrow\P^2$ is a family of genus
two curves. Assume this family has mild degenerations and let
$\pi:X\rightarrow\P^2$ be the compactified relative Jacobian of
$Y/\P^2$. The previous lemma says that $W$ must be isomorphic to
$$R^1f_*\O_Y\cong R^1\pi_*\O_X$$
up to tensoring with a line bundle. We can identify this line bundle
in terms of the branch locus of the double cover. This is essentially
Proposition~5 of Markushevich~\cite{markushevich95}, though we avoid
local coordinates in our proof.

\begin{lem}
\label{normalize}
In the above situation
$$W\cong {\O}(d)\otimes \Lambda^2W^*\otimes R^1f_*\O_Y.$$
\end{lem}

\begin{prf}
Denote the double cover by $g:Y\rightarrow\P(W)$. We have a
commutative diagram
$$\begin{array}{crc}
	Y & \stackrel{g}{\rightarrow} & {\P}(W) \\
	 & f\searrow & \downarrow h \\
	 & & {\P}^2. \\
	 \end{array}$$
Then
$$R^1f_*\O_Y\cong R^1(h\circ g)_*\O_Y\cong R^1h_*\circ R^0g_*\O_Y$$
since $g$ is a finite map and therefore $R^1g_*\O_Y$ vanishes. Because
of the choice of branch divisor
$$R^0g_*\O_Y\cong\O_{\P(W)}\oplus\left(\O_{\P(W)}(-3)\otimes h^*\O(-d)\right),$$
and thus
\begin{eqnarray*}
R^1f_*{\O}_Y & \cong & R^1h_*{\O}_{{\P}(W)}\oplus R^1h_*\left({\O}_{{\P}(W)}(-3)\otimes h^*{\O}(-d)\right) \\
 & \cong & {\O}(-d)\otimes R^1h_*{\O}_{{\P}(W)}(-3) \\
 & \cong & {\O}(-d)\otimes \left(R^0h_*\omega_h(3)\right)^*
 \end{eqnarray*}
where on the second line we have used the fact that $R^1h_*\O_{\P(W)}$
vanishes (the fibres of $h:\P(W)\rightarrow\P^2$ are rational curves) and
on the third line we have used relative Serre duality. The relative dualizing sheaf $\omega_h$ of $h$ is given by
$$\omega_h\cong \O_{\P(W)}(-2)\otimes h^*\Lambda^2W^*.$$
Inserting this in the above equation gives
\begin{eqnarray*}
R^1f_*{\O}_Y & \cong & {\O}(-d)\otimes \left(R^0h_*({\O}_{{\P}(W)}(1)\otimes h^*\Lambda^2W^*)\right)^* \\
 & \cong & {\O}(-d)\otimes \Lambda^2W\otimes \left(R^0h_*{\O}_{{\P}(W)}(1)\right)^* \\
 & \cong & {\O}(-d)\otimes \Lambda^2W\otimes W
 \end{eqnarray*}
and the statement of the lemma follows.
\end{prf}

We have already established which rank two vector bundles can arise as
$R^1\pi_*\O_X$ for a four-fold with trivial canonical bundle $X$
fibred over $\P^2$ by abelian surfaces. Recall that $R^1\pi_*\O_X$
must be one of the following:
\begin{enumerate}
\item $\O_{\P^2}\oplus\O_{\P^2}(-3)$ if the finite cover
  $\tilde{X}$ of $X$ is either an abelian four-fold, a product
  $A\times S$ of abelian and K3 surfaces, or a product $E\times Y$ of
  an elliptic curve and a Calabi-Yau three-fold,
\item $\O_{\P^2}(-1)\oplus\O_{\P^2}(-2)$ or
  $\O_{\P^2}(-2)\oplus\O_{\P^2}(-2)$ if $\tilde{X}$ is a Calabi-Yau
  four-fold,
\item $\Omega^1_{\P^2}$ if $\tilde{X}=X$ is an irreducible
  holomorphic symplectic four-fold.
\end{enumerate}
If $X$ is the compactified relative Jacobian of a family of curves
with mild degenerations, then the generic singular fibre of
$X\rightarrow\P^2$ will be the compactified Jacobian of a curve with a
single node. In this case, both of the conditions in
Lemma~\ref{first_chern} will be satisfied, implying $c_1(R^1\pi_*\O_X)=-3$. Thus
$\O_{\P^2}(-2)\oplus\O_{\P^2}(-2)$ is impossible; we are left
with the three other possibilities.

Returning to the construction of the family of curves, $W$ must be
isomorphic to $R^1f_*\O_Y$ up to tensoring with a line bundle. We can
normalize the first Chern class of $W$ by assuming this line bundle is
trivial, i.e., we set
$$W\cong R^1f_*\O_Y\cong R^1\pi_*\O_X.$$
The isomorphism
$$W\cong {\O}(d)\otimes \Lambda^2W^*\otimes R^1f_*\O_Y$$
of Lemma~\ref{normalize} then yields
$$d=c_1(W)=c_1(R^1f_*\O_Y)=c_1(R^1\pi_*\O_X)=-3,$$
and hence the branch divisor is the zero locus of a generic section of
$$\O_{\P(W)}(6)\otimes h^*\O(-6)$$
on $\P(W)$. Now
\begin{eqnarray*}
{\H}^0\left({\P}(W),{\O}_{{\P}(W)}(6)\otimes h^*{\O}(-6)\right) & \cong & {\H}^0\left({\P}^2,R^0h_*\left({\O}_{{\P}(W)}(6)\otimes h^*{\O}(-6)\right)\right) \\
 & \cong & {\H}^0\left({\P^2},{\O}(-6)\otimes R^0h_*{\O}_{{\P}(W)}(6)\right) \\
 & \cong & {\H}^0\left({\P^2},{\O}(-6)\otimes \mathrm{Sym}^6W^*\right).
\end{eqnarray*}

We consider the three possibilities for $W\cong R^1f_*\O_Y\cong R^1\pi_*\O_X$.

\begin{prp}
\label{no0-3}
The case $W\cong\O_{\P^2}\oplus\O_{\P^2}(-3)$ is impossible. In
particular, if the finite cover $\tilde{X}$ of $X$ is either an
abelian four-fold, a product $A\times S$ of abelian and K3 surfaces,
or a product $E\times Y$ of an elliptic curve and a Calabi-Yau
three-fold then $X$ cannot be a fibration by Jacobians over $\P^2$.
\end{prp}

\begin{prf}
If $W\cong\O_{\P^2}\oplus\O_{\P^2}(-3)$ then the branch divisor of
$Y\rightarrow\P(W)$ corresponds to a section in
$$\H^0(\P^2,{\O}(-6)\otimes
\mathrm{Sym}^6W^*)\cong\H^0(\P^2,{\O}(-6)\oplus {\O}(-3)\oplus\ldots\oplus {\O}(12)).$$
If we write this section as $(s_0,s_1,\ldots,s_6)$ then $s_0$
and $s_1$ must vanish, since they are sections of line bundles of
negative degree. Let us consider what this means.

A fibre $Y_t$ of $Y\rightarrow\P^2$ is a genus two curve, which is a
double cover of $\P(W_t)$ branched over six points. Each point
corresponds to a line in $W_t^*$, and therefore the six points define
a line in $\mathrm{Sym}^6W_t^*$. This is the same as the line in
$\mathrm{Sym}^6W_t^*$ given by evaluating the section
$(s_0,s_1,\ldots,s_6)$ at $t$, since the branch locus is defined by
this section. We recover the six points from $(s_0,s_1,\ldots,s_6)$ by
taking the roots of the polynomial
$$s_6z^6+s_5z^5+\ldots+s_1z+s_0.$$
If $s_0$ and $s_1$ vanish at $t$, then $z=0$ is a repeated root. This
means that two of the branch points coincide and the genus two
curve $Y_t$ is singular. If $s_0$ and $s_1$ vanish identically, then
all curves $Y_t$ in the family $Y\rightarrow\P^2$ are singular. This
contradicts the mild degenerations assumption.
\end{prf}

\begin{prp}
\label{CY4-fold}
The case $W\cong\O_{\P^2}(-1)\oplus\O_{\P^2}(-2)$ produces a
Calabi-Yau four-fold $X$. This example belongs to a 75-dimensional
family.
\end{prp}

\begin{rmk}
Markushevich described this four-fold in Example~1
of~\cite{markushevich95}. Our results enable us to identify it as a
Calabi-Yau four-fold.
\end{rmk}

\begin{prf}
If $W\cong\O_{\P^2}(-1)\oplus\O_{\P^2}(-2)$ then the branch divisor of
$Y\rightarrow\P(W)$ corresponds to a section in
$$\H^0(\P^2,{\O}(-6)\otimes
\mathrm{Sym}^6W^*)\cong\H^0(\P^2,{\O}\oplus {\O}(1)\oplus\ldots\oplus {\O}(6)).$$
If this section is chosen generically then the resulting family of
genus two curves $Y\rightarrow\P^2$ will have mild degenerations, and
hence by Theorem~2 of Markushevich~\cite{markushevich95} its
compactified relative Jacobian $X$ will be a smooth four-fold whose
canonical bundle $K_X$ is trivial when restricted to the fibres of
$\pi:X\rightarrow\P^2$. In particular, $K_X$ is the pull-back
$\pi^*\O_{\P^2}(m)$ of a line bundle on $\P^2$. We can proceed as in
the example of Section~7: substituting
$$R^1\pi_*\O_X\cong W\cong \O_{\P^2}(-1)\oplus\O_{\P^2}(-2)$$
into the Leray spectral sequence gives
$$h^k(\O_X)=\left\{ \begin{array}{ll}
  1 & \mbox{if }k=0\mbox{ or }4, \\
  0 & \mbox{otherwise.} \\
\end{array}\right.$$
Therefore $m=0$ and $K_X$ is trivial. A four-fold $X$ with trivial
canonical bundle and these Hodge numbers must be a Calabi-Yau
four-fold or the quotient of a product $S_1\times S_2$ of K3 surfaces
by an involution, but by Proposition~\ref{impossible2} the latter
cannot be fibred by abelian surfaces over $\P^2$.

The space of sections above has dimension
$$1+3+6+10+15+21+28=84.$$
Projectivizing and quotienting by the automorphism group
$\mathrm{PGL}(3,\C)$ of $\P^2$ leaves
$$84-1-8=75$$
parameters.
\end{prf}

The final result of this section is really Theorem~5 of
Markushevich~\cite{markushevich95}; we include a slight simplification
of his proof, which avoids calculations in local coordinates.

\begin{prp}
\label{hilb2K3}
The case $W\cong\Omega^1_{\P^2}$ produces an irreducible holomorphic
symplectic four-fold $X$. Moreover, $X$ is isomorphic to the
Beauville-Mukai system, i.e., the family of genus two curves is a
complete linear system of curves on a K3 surface $S$. In particular,
$X$ is deformation equivalent to $\mathrm{Hilb}^2S$. This example
belongs to a 19-dimensional family.
\end{prp}

\begin{prf}
If $W\cong\Omega^1_{\P^2}$ then the branch divisor of
$Y\rightarrow\P(W)$ corresponds to a section in
$$\H^0(\P^2,{\O}(-6)\otimes
\mathrm{Sym}^6W^*)\cong\H^0(\P^2,{\O}(-6)\otimes \mathrm{Sym}^6T).$$
As in the previous proposition, a generic choice of section produces a
family of curves with mild degenerations whose compactified relative
Jacobian $X$ is a smooth four-fold whose canonical bundle $K_X$ is the
pull-back of a line bundle on $\P^2$. Since 
$$R^1\pi_*\O_X\cong W\cong\Omega^1_{\P^2},$$
the Leray spectral sequence gives
$$h^k(\O_X)=\left\{ \begin{array}{ll}
  1 & \mbox{if }k=0,2,\mbox{ or }4, \\
  0 & \mbox{otherwise.} \\
\end{array}\right.$$
As before, we can argue that $K_X$ is trivial. Then the holonomy
classification (Proposition~\ref{berger}) implies that $X$ is an
irreducible holomorphic symplectic four-fold. In particular, the
generator $\sigma$ of the one-dimensional space
$$\H^0(\Omega^2_X)\cong \H^2(\O_X)$$
is a holomorphic symplectic form, and $X$ is simply connected.

Next we calculate the dimension of the space of sections
$$\H^0(\P^2,{\O}(-6)\otimes \mathrm{Sym}^6T).$$
Writing $\P^2$ as a symmetric space $\mathrm{GL}(3,\C)/P$, the
homogeneous bundle
$$\O(-6)\otimes\mathrm{Sym}^6T$$
corresponds to the $P$-representation $(0|0,6)$ (see Eastwood and
Sawon~\cite{es02} for notation). The Borel-Weil Theorem then implies
that
$$\H^0(\P^2,{\O}(-6)\otimes \mathrm{Sym}^6T)$$
is isomorphic to the $\mathrm{GL}(3,\C)$-representation
$$(0,0,6)\cong\mathrm{Sym}^6\C^3.$$
This has dimension 28. After subtracting one (rescaling of the
section) and eight (quotienting by the action of
$\mathrm{PGL}(3,\C)$), we are left with 19 parameters for this family
of examples.


Finally, we relate this example to the Beauville-Mukai system. Observe
that
$$\P(W)\cong\P(\Omega^1_{\P^2})$$
is the incidence subvariety in $\P^2\times (\P^2)^*$: a point $q$ in
$(\P^2)^*$ corresponds to a line in $\P^2$, and $(p,q)\in\P(W)$ if and
only if $p$ lies on this line. We previously used the isomorphism
$${\H}^0\left({\P}(W),{\O}_{{\P}(W)}(6)\otimes
h^*{\O}(-6)\right)\cong{\H}^0\left({\P^2},{\O}(-6)\otimes
\mathrm{Sym}^6W^*\right).$$
However, in this example there is also a projection $\P(W)\rightarrow
(\P^2)^*$. We claim that the line bundle
$${\O}_{{\P}(W)}(6)\otimes h^*{\O}(-6)$$ 
on $\P(W)$ is the pull-back of the line bundle $\O_{(\P^2)^*}(6)$ on
$(\P^2)^*$. This follows directly from the fact that
${\O}_{{\P}(W)}(1)$ is the restriction to
$$\P(W)\subset \P^2\times (\P^2)^*$$
of the line bundle $\O_{\P^2}(1)\boxtimes\O_{(\P^2)^*}(1)$ on
$\P^2\times (\P^2)^*$. Therefore
$${\H}^0\left({\P}(W),{\O}_{{\P}(W)}(6)\otimes
h^*{\O}(-6)\right)\cong{\H}^0\left((\P^2)^*,{\O}_{(\P^2)^*}(6)\right)$$
and the branch divisor is the pull-back of a sextic curve in $(\P^2)^*$.

Let $S$ be the double cover of $(\P^2)^*$ branched over the
sextic. For a generic choice of sextic, $S$ will be a (smooth) K3
surface; $Y$ will be a $\P^1$-bundle over $S$. We have the following
commutative diagram.
$$\begin{array}{ccccc}
 & & Y & & \\
 & \swarrow & \downarrow & & \\
S & & {\P}(W) & & \\
\downarrow & \swarrow & & \searrow & \\
({\P}^2)^* & & & & {\P}^2 \\
\end{array}$$
If $t\in\P^2$ then $Y_t$ is a genus two curve which is a double cover
of $\P(W_t)$. The line $\P(W_t)$ projects isomorphically to a line in
$(\P^2)^*$, whose preimage in $S$ will be a curve isomorphic to
$Y_t$. Therefore the family of curves $Y\rightarrow\P^2$ is
(isomorphic to) a complete linear system of curves in the K3 surface
$S$. This completes the proof.
\end{prf}

\begin{rmk}
In Section~6 of~\cite{markushevich96}, Markushevich described a method
of calculating
$$\H^0(B,K_B^2\otimes \mathrm{Sym}^6T_B)$$
for a toric variety $B$. This was important because he used this
method to show that $\P^2$ is the only rational surface that can arise
as the base of a fibration on an irreducible holomorphic symplectic
four-fold. However, later results of Matsushita~\cite{matsushita99}
allow us to conclude more directly that the base must be $\P^2$, so we
instead use the Borel-Weil Theorem to calculate the above space of
sections.
\end{rmk}

\begin{rmk}
If $S$ is a double cover of the projective plane branched over a
sextic, then it has a polarization $H$ with $H^2=2$. The moduli space
of polarized K3 surfaces with $H^2=2$ is 19-dimensional. For each of
these we get a Beauville-Mukai system
$$X:=\overline{\mathrm{Jac}}^0(\mathcal{C}/\P^2).$$
Note that the dimension of this family of examples agrees with that
computed in Proposition~\ref{hilb2K3}.
\end{rmk}

\begin{flushleft}
Department of Mathematics\hfill sawon@math.colostate.edu\\
Colorado State University\\
Fort Collins CO 80523-1874\\
USA\\
\end{flushleft}

\end{document}